\documentclass [11pt,reqno]{amsart}
\usepackage {amsmath,amssymb,verbatim,geometry}
\usepackage[all]{xy}
\newif\ifpdf
\ifpdf
\usepackage[pdftex]{hyperref}
\else
\fi

\numberwithin{equation}{section}       

\newtheorem{prop} {Proposition} [section]
\newtheorem{thm}[prop] {Theorem} 
\newtheorem{quest} {Question} 
\newtheorem{defi}[prop] {Definition}

\newtheorem{prop-def}[prop]{Proposition-Definition}

\newtheorem{exa}[prop]{Example}

\theoremstyle{remark}

\newtheorem*{ackn}{Acknowledgements}

\newcommand{\C}{{\mathbb{C}}}
\newcommand{\B}{{\mathbb{B}}}

\newcommand{\PP}{{\mathbb{P}}}

\newcommand{\R}{{\mathbb{R}}}

\newcommand{\cH}{{\mathcal{H}}}

\newcommand{\cO}{{\mathcal{O}}}

\renewcommand{\a}{\alpha}

\newcommand{\e}{\varepsilon}
\newcommand{\om}{\omega}
\newcommand{\f}{\varphi}
\newcommand{\p}{\psi}

\newcommand{\D}{\Delta}

\newcommand{\MA}{\mathrm{MA}}

\newcommand{\Ric}{\mathrm{Ric}}

\newcommand{\psh}{{\mathrm{PSH}}}

\title{Open problems in pluripotential theory}

\date{\today}

\author{ S.Dinew*, V.Guedj**
    \and
    A. Zeriahi**}

 \address{Jagiellonian University 30-348 Krakow, Lojasiewicza 6, Poland}
\email{slawomir.dinew@im.uj.edu.pl}

 \address{Institut Universitaire de France \& Institut Math\'ematiques de Toulouse, Toulouse,
France}
\email{vincent.guedj@math.univ-toulouse.fr}

\address{Institut Math\'ematiques de Toulouse, Universit{\'e} Paul Sabatier, Toulouse
France}
\email{ahmed.zeriahi@math.univ-toulouse.fr}

\subjclass[2010]{Primary: 32W20, seconday: 32Q20 \\
* This author is partially supported by grant NCN 2013/11/D/ST1/02599\\
**These authors are partially supported by the french A.N.R. project GRACK}

\begin{document}

\begin{abstract}  
We propose a list of open problems in pluripotential theory partially motivated by their applications  to complex differential geometry. The list includes both local questions as well as issues related to the compact complex manifold setting.
\end{abstract} 

\maketitle

\setcounter{tocdepth}{1}

\tableofcontents

\newpage

\section*{Introduction}

Pluripotential theory is the several complex variables analogue of classical potential theory in the complex plane (or in Riemann surfaces). While the latter is a linear theory associated to the Laplacian
of a K\"ahler metric, the former is highly non-linear and associated to 
complex Monge-Amp\`ere operators
$$
\f \mapsto MA(\f):=c \, (\omega+dd^c \f)^n,
$$
where $\omega$ is a K\"ahler metric (or merely a semi-positive and big form), $n$ is the complex dimension
of the ambient manifold
and $c>0$ a normalizing constant.

When $\f$ is smooth and the ambient manifold is an open subset of 
$\C^n$ equipped with the euclidean K\"ahler metric 
$\omega=dd^c ||z||^2=\sum_{j=1}^n \frac{i}{\pi} dz_j \wedge d \overline{z}_j$, a change of unknown
function $\f \mapsto \f+||z||^2$ reduces the previous expression to
$$
MA(\f)=c' \, \det \left( \frac{\partial^2 \f}{\partial z_i \partial \overline{z}_j} \right) \, dV(z),
$$
where $dV$ denotes the Lebesgue volume form and $c'>0$ is  a normalizing constant.

Complex Monge-Amp\`ere equations 
$MA(\f)=F(\f) \mu$
have been one of the most powerful tools in K\"ahler geometry since Aubin and Yau's classical works 
\cite{Aub78,Yau78}, culminating in Yau's solution to the Calabi conjecture. A notable application is the construction of K\"ahler-Einstein metrics on compact K\"ahler manifolds.

In recent years, following Tsuji's pioneering work \cite{Tsu88}, degenerate complex Monge-Amp\`ere equations have been intensively studied by many authors. In relation to the Minimal Model Program, they led to the construction of singular K\"ahler-Einstein metrics
\cite{EGZ09,BBGZ13,BBEGZ} or, more generally, of canonical volume forms on compact K\"ahler manifolds with nonnegative Kodaira dimension [ST12].

\smallskip

Making sense of and constructing singular K\"ahler-Einstein metrics on these mildly singular varieties has required advanced tools in the study of degenerate complex Monge-Amp\`ere equations and led
to many open problems.

 We have been asked by several colleagues to propose a list of such open questions.
 There is no pretention to be exhaustive, 
the selection we propose reflects our own taste and limitations. 
The questions are of various nature, some would lead to breakthrough and require new ideas,
others are probably reasonably easy.

Although there are many interesting subjects that we haven't  covered, we hope there is nevertheless 
enough interesting questions for the curious reader to work on.

\subsection*{Structure of the article}
The article is organized as follows:
\begin{itemize}
\item In the first section we list open questions related to the domain of definition of the complex Monge-Amp\`ere operator and quasi-plurisubharmonic Green functions.
\item    We collect in the second section problems related to the description of the range of the complex Monge-Amp\`ere operator acting on various classes of (quasi-)psh functions.
\item  The third section focuses on regularity issues, in the local and  the compact settings.
\item The fourth section is devoted to pluripotential theory on non K\"ahler  manifolds.
\item We mention in the fifth section some problems about the extension of quasi-psh functions,
as well as how to control their singularities.
\item The sixth and last section is devoted to various open questions that arise in the viscosity approach,
for degenerate elliptic and parabolic Monge-Amp\`ere equations.
\end{itemize}

\begin{ackn} 
The authors would like to express their gratitude to Z.Blocki, D.Coman, U.Cegrell C.Favre, S.Ko\l{}odziej, M.P\u aun and A.Rahskovskii for several useful conversations.
\end{ackn}

\section{Green functions}

We let $MA(\f)$ denote the complex Monge-Amp\`ere measure of a (quasi-)plurisubharmonic function.
When $\f$ is  smooth  and plurisubharmonic in a domain $\Omega \subset \C^n$, then
$$
MA(\f)=(dd^c \f)^n=c \, \det \left( \frac{\partial^2 \f}{\partial z_i \partial \overline{z}_j} \right) \, dV(z),
$$
where $dV$ denotes the Lebesgue volume form and $c>0$ is  a normalizing constant.

When $\f$ is smooth and $\omega$-plurisubharmonic on a compact K\"ahler manifold $(X,\omega)$, then
$$
MA(\f)=c \, (\omega+dd^c \f)^n,
$$
where $c=(\int_X \omega^n)^{-1}$ is chosen so that $MA(\f)$ is a probability measure.

\subsection{Domains of definition}

\subsubsection{The local context}

When $n=1$, the complex Monge-Amp\`ere operator reduces to the linear  Laplace operator $\Delta$, whose action on all (p)sh functions is well defined in the sense of distributions. In higher dimensions one cannot usually define the product of two distributions, so the definition of the complex Monge-Amp\`ere operator is more involved.

\begin{defi}
Let $\f$ be a plurisubharmonic function in some domain $\Omega \subset \C^n$. We say that
the complex Monge-Amp\`ere measure $(dd^c \f)^n$ is well defined if there exists a Radon measure
$\mu$ such that for any sequence 
$(\f_j)$ of smooth plurisubharmonic functions decreasing to $\f$ in some relatively compact
subdomain $V \subset \subset \Omega$, the Monge-Amp\`ere measures
$(dd^c \f_j)^n$ converge weakly to $\mu$ in $V$.
If such is the case, one sets $(dd^c \f)^n:=\mu$.
\end{defi}

Bedford and Taylor have shown that the complex Monge-Amp\`ere operator $\f \mapsto (dd^c \f)^n$ is well defined for all {\it locally bounded}
psh functions \cite{BT82}, however it can not be  defined for all
psh functions when $n \geq 2$ (see \cite{Kis83}), motivating the following questions \cite[p69]{Bed93}:

 \begin{quest} [Bedford]
 Given $\f$ an arbitrary plurisubharmonic function in some domain $\Omega \subset \C^n$, what can be said about the set of cluster points
 of $(dd^c \f_j)^n$, where $\f_j$ are bounded plurisubharmonic functions which decrease to $\f$ ?
\end{quest}

By definition this set is reduced to a single Radon measure when $\f$ belongs to the domain of definition 
$DMA(\Omega)$ of the complex Monge-Amp\`ere operator. One further difficulty
is that this operator is not continuous for the $L^1_{loc}$-topology \cite{Ceg83,Lel83}:

 \begin{quest} [Bedford]
 Given $\f$ a bounded plurisubharmonic function in some domain $\Omega \subset \C^n$, what can be said about the set of cluster points of
 $(dd^c \f_j)^n$, where $\f_j$ are bounded plurisubharmonic functions which converge to $\f$ in $L_{loc}^1$ ?
\end{quest}

When  $n=2$, a psh function $\f$ belongs to the domain of definition $DMA(\Omega)$,
$\Omega$ a bounded domain,
if and only if its gradient 
$\nabla \f$ belongs to $L^2_{loc}(\Omega)$, 
$$
DMA(\Omega)=PSH(\Omega) \cap W^{1,2}_{loc}.
$$
Indeed if such is the case, a simple integration by parts shows that $\f$ is locally integrable with respect to the trace measure of the positive current $dd^c \f$,
hence the current $\f \, dd^c \f$ is well defined and so is
$$
(dd^c \f)^2:=dd^c (\f dd^c \f).
$$
This has been observed by Bedford and Taylor in \cite{BT78}.
To show the converse, Cegrell provides in \cite{Ceg02} explicit computations showing that 
$-\sqrt{-\log|z_2|}$ does not belong to $DMA(\B)$, as it can be approximated by a decreasing sequence of bounded psh functions whose Monge-Amp\`ere masses blow up.
This is extended by Blocki  in \cite{Blo04} who shows that  {\it any} psh function 
$\f \notin W^{1,2}_{loc}(\Omega)$ is the decreasing limit of smooth psh functions whose Monge-Amp\`ere masses blow up.
A simpler argument has been further provided by Cegrell in \cite{Ceg07}.

Let us recall that psh functions have gradient in $L^q_{loc}$ for all $q<2$. The function 
$\f(z)=\log|z_1|$ does not have gradient in $L_{loc}^2$ while any locally bounded psh function $\f$ 
does: indeed we can assume $\f \geq 0$ and observe that 
$$
dd^c \f^2 =2 \f dd^c \f +2 d\f \wedge d^c \f,
$$ 
hence the $L^2$-norm of 
$\nabla \f$ is bounded from above by the mass of the positive current $dd^c \f^2$.

In higher dimension $n \geq 3$, the domain of definition  $DMA(\Omega)$   has been fully characterized by Cegrell \cite{Ceg04} and Blocki \cite{Blo06}, but the characterization is more involved.
We only stress, for further use, the following important subclass: plurisubharmonic functions with 
{\it compact singularities} belong to $DMA(\Omega)$. We say here that a function $\f$ has compact singularities if it is locally bounded near the boundary of $\Omega$. 

\subsubsection{Maximal plurisubharmonic functions}

\begin{defi}
A plurisubharmonic function $\f$ is called maximal in a domain $\Omega \subset \C^n$
if for any other plurisubharmonic function $v$ in $\Omega$ such that $v \leq \f$ near $\partial \Omega$,
then $v \leq \f$ in $\Omega$.
\end{defi}

When $n=1$, maximal plurisubharmonic function are precisely  harmonic ones, in particular they are smooth, but in higher dimension they may be quite irregular:  a plurisubharmonic function that only depends on
$n-1$ variables is maximal. When $\f \in DMA(\Omega)$ then
$$
\f \text{ is maximal } \Longleftrightarrow MA(\f)=0,
$$
as explained by Blocki in \cite{Blo04}. In particular maximality is a local notion in $DMA(\Omega)$.

 \begin{quest} [Blocki]
Is maximality a local notion for arbitrary plurisubharmonic functions?
\end{quest}

If there exists a holomorphic foliation by complex discs $\Delta_t$ along which a  plurisubharmonic function
$\f$ is harmonic, then $\f$ is maximal (see e.g. \cite{DG12}).

Conversely assume $n=2$ (for simplicity): if a maximal plurisubharmonic function $\f$ is ${\mathcal C}^3$-smooth, it follows from Frobenius theorem that
the support of the differential form $dd^c \f$ is foliated by holomorphic discs $\Delta_t$ such that $\f_{|\Delta_t}$
is harmonic (see \cite{BK77}). Dujardin has provided in \cite{Duj10} several examples of maximal 
plurisubharmonic functions that are ${\mathcal C}^{1,\alpha}$-smooth and such that 
there is no such holomorphic disc.

 \begin{quest} [Dujardin]
 Assume  $0 \leq \alpha <1$ and $n=2$.
Assume $\f$ is ${\mathcal C}^{2,\alpha}$-smooth and maximal. Is the suppport of $dd^c \f$ foliated (or laminated)
by holomorphic disc along which $\f$ is harmonic ?
\end{quest}

\subsubsection{Compact K\"ahler manifolds}

The situation turns out to be quite different on compact K\"ahler manifold.
Let $(X,\omega)$ be a compact K\"ahler manifold and let $PSH(X,\omega)$ denote the set of all 
$\omega$-psh functions. By mimicking the definition in the local case, we propose:

\begin{defi}
Fix $\f \in PSH(X,\omega)$. We say that
the complex Monge-Amp\`ere measure $(\omega+dd^c \f)^n$ is well defined and write 
$\f \in DMA(X,\omega)$, if 
there exists a Radon measure $\mu$ such that for any sequence 
$(\f_j)$ of smooth $\omega$-plurisubharmonic functions decreasing to $\f$  
on $X$, the Monge-Amp\`ere measures
$(\omega+dd^c \f_j)^n$ converge weakly to $\mu$.

If such is the case, one sets $(\omega+dd^c \f)^n:=\mu$.
\end{defi}

The existence of smooth $\omega$-psh approximants $\f_j \in PSH(X,\omega) \cap {\mathcal C}^{\infty}(X)$
is not an obvious fact, see \cite{Dem92,BK07}. One can replace in the definition {\it smooth}  by
{\it bounded} approximants by using Bedford-Taylor's result: bounded $\omega$-psh approximants are easy to construct, for example $\f_j=\max(\f,-j) \in PSH(X,\omega) \cap L^{\infty}(X)$
would do.

\begin{defi}
We say that $\f$ belongs to $DMA_{loc}(X,\omega)$ if for any 
  covering of $X$ by bounded strictly pseudonvex domains $U_{\alpha}$, and any local potential
$\rho_{\alpha}$ of $\omega$ in $U_{\alpha}$, the function
$\rho_{\alpha}+\f$ belongs to $DMA(U_{\alpha})$.
\end{defi}

The (global) domain of definition $DMA(X,\omega)$ of the complex Monge-Amp\`ere operator
$\f \mapsto (\omega+dd^c \f)^n$ is different than its local counterpart
$DMA_{loc}(X,\omega)$.
The construction of Blocki \cite{Blo04} breaks down here as the total mass of
the  Monge-Amp\`ere measures of any approximating sequence is fixed, equal to the volume 
$Vol_\om(X)=\int_X \omega^n$.
It turns out indeed that $DMA(X,\omega)$ contains functions
which do not have gradient in $L^2$.

 \begin{quest} [Coman-Guedj-Zeriahi]
Characterize $DMA(X,\omega)$.
\end{quest}

Again $DMA(X,\omega)$ is a strict subclass of $PSH(X,\omega)$, Lelong numbers provide for instance an obstruction to being in $DMA(X,\omega)$ (see \cite{CGZ08}).
An important subclass ${\mathcal E}(X,\omega) \subset DMA(X,\omega)$ of functions with full 
Monge-Amp\`ere mass has been introduced in \cite{GZ07} (see section \ref{sec:fullMA}).
Functions in ${\mathcal E}(X,\omega)$ have zero Lelong numbers. Thus $DMA(X,\omega)$ contains
the two distinct subsets ${\mathcal E}(X,\omega)$ and $DMA_{loc}(X,\omega)$. 
Apart from this not much is known:

 \begin{quest} [Coman-Guedj-Zeriahi]
Is $DMA(X,\omega)$ convex ? 
If $\f \in DMA(X,\omega)$ and $\p \in PSH(X,\omega)$, does it imply that $\max(\f,\p) \in DMA(X,\omega)$ ?
\end{quest}

We refer the reader to \cite{CGZ08} for more open questions, notably a tentative candidate for 
describing $DMA(X,\omega)$ in complex dimension $n=2$.

\subsection{The Dirac mass question} \label{sec:Dirac}

\text{ }

Let $\f$ be a plurisubharmonic function in the unit ball $\B$ of $\C^n$. Assume $\f$ has a well-defined 
complex Monge-Amp\`ere measure $(dd^c \f)^n$, we could for instance 
assume that $\f$ is locally bounded  outside the origin.

\begin{quest} [Guedj, Rashkovskii] \label{q:GR}
Assume that $(dd^c \f)^n $ has a Dirac mass at the origin. Does it imply that 
$\f$ has a positive Lelong number at $0 \in \C^n$ ?
\end{quest}

The reverse implication is well-known to hold (see \cite{Dem93}). Note however
that 
$$
\f(z,w)=\max( \e \log |z|, \e^{-1} \log |w|)
$$
yields an example of a psh function such that $(dd^c \f)^2$  has Dirac mass $1$ at the origin, but 
an arbitrarily small Lelong number $\nu(\f,0)=\e$.

We list below 
some partial information and a few situations where the answer is known:
\begin{itemize}
\item one necessarily has $\f(0)=-\infty$ (an observation due to Cegrell, see \cite{Ceg02});
\item  if $\f$ is {\it toric}, i.e. $\f(e^{i\theta_1} z_1,\ldots,e^{i\theta_n} z_n)=\f(z_1,\ldots,z_n)$
for all $\theta=(\theta_1,\ldots,\theta_n)$ and $z=(z_1,\ldots,z_n)$, then the answer is positive (see \cite{Ra01});
\item  if $e^{\f}$ is H\"older continuous, in particular if $\f$ has analytic singularities, then the answer is positive (see \cite{FJ07,Ra13});
\item  if $\f$ is bounded from below by $\gamma \log \| z\|-C$, the answer is positive \cite{Wik05};

\item when $\f$ is invariant under a rational mapping, the answer is positive 
(see \cite{G10}).
\end{itemize}

Demailly has proposed \cite{Dem14} a stronger question. Consider the Bergman spaces
$$
\cH(j \f):=\left\{ f \in \cO(\B) \, | \, \int_\B |f|^2 e^{-2j\f} dV <+\infty \right\}
$$
of holomorphic functions that are $L^2$ in the ball with respect to the measure $e^{-2j\f} dV$, where
$dV$ denotes the Lebesgue measure. Let $(f_{k,j})$ denote an orthonormal basis of $\cH(j \f)$ and set
$$
\f_j:=\frac{1}{2j} \log \left[ \sum_{k=0}^{+\infty} |f_{k,j}|^2 \right].
$$

As observed by Demailly in \cite{Dem93}, it follows from the Ohsawa-Takegoshi theorem that 
the $\f_j$'s approximate $\f$ in a very accurate manner, hence motivating the following:

 \begin{quest} [Demailly]
Does one have $(dd^c \f_j)^n(0) \longrightarrow (dd^c \f)^n(0)$  ?
\end{quest}

Since the $\f_j$'s have analytic singularities and $\nu(\f_j,0) \rightarrow \nu(\f,0)$, a positive answer to Demailly's question would also imply a positive answer to Question \ref{q:GR}.
The interested reader will get more information on these questions in Rahskovskii's papers (see
\cite{Ra06,Ra13}).

\subsection{Quasi-plurisubharmonic Green functions}

Cegrell and Wiklund have made some contributions to Question \ref{q:GR}
showing in particular that it suffices to answer it for functions whose Monge-Amp\`ere measure
is concentrated at the origin (see \cite{Wik05}).
Such functions are called {\it Green functions} whenever they have zero boundary values.

Classical Green functions are fundamental solutions to the Laplace equations. As such they play a crucial role in  linear potential theory. 
Appropriate notions of plurisubharmonic Green functions have been defined and studied in domains of $\C^n$, we refer the reader to \cite{Kis00,PSS12} for a survey of their properties and applications, see also section 5.6.

One can similarly study $\omega$-psh Green functions on a compact K\"ahler manifold 
$(X,\omega)$: for $a \in X$, these are the functions $g_a \in DMA(X,\omega)$ such that
$$
MA(g_a)=V^{-1}(\omega+dd^c g_a)^n=\delta_a
$$
is the Dirac mass at point $a$. 
Such functions have been used in complex dynamics to construct invariant measures of maximal entropy
(see \cite{CG04,G05,G10}). 

The article \cite{CG09} studies quasi-plurisubharmonic Green functions for themselves, exhibits
some examples and several questions, for instance:

 \begin{quest} [Coman-Guedj]
Assume $(X,\omega)$ is the complex projective space $\C\PP^n$ equipped with the Fubini-Study K\"ahler form. Can one find examples of a $\omega$-psh function $G_a \in DMA(X,\omega)$, whose Monge-Amp\`ere measure is concentrated at some point $a \in \C\PP^n$ and having an irrational Lelong number at $a$ ?
\end{quest}

The answer is negative in dimension $1$, as the Lelong number coincides in this case with the Dirac mass.
The situation is more subtle in higher dimension, as explained in section \ref{sec:Dirac}.
The reader will find in \cite{CG09} various constructions of such  quasi-psh Green functions with given rational Lelong number at $a$.

\section{Range of the complex Monge-Amp\`ere operator}  

In this section we try and describe the range of the complex Monge-Amp\`ere operator on various classes
of  ($\omega$-)plurisubharmonic functions.

\subsection{The non-pluripolar product}  \label{sec:fullMA}

Let $(X,\omega)$ be a compact K\"ahler manifold. Given $\f \in PSH(X,\omega)$ we consider the bounded approximants
$$
\f_j:=\max(\f,-j) \in PSH(X,\omega) \cap L^{\infty}(X).
$$
It has been realized in \cite{GZ07} that the truncated Monge-Amp\`ere measures
$$
1_{\{ \f>-j\}} (\omega+dd^c \f_j)^n
$$
form an increasing sequence of Borel measures with  mass uniformly bounded by $V=\int_X \omega^n$.
A similar local construction had been previously proposed by Bedford and Taylor  \cite{BT87}.

\begin{defi}
Set $T=\omega+dd^c \f \geq 0$. The measure 
$$
\langle T^n \rangle=\lim \nearrow  1_{\{ \f>-j\}} (\omega+dd^c \f_j)^n
$$
is called the non-pluripolar measure product of $T^n$.
\end{defi}

It has been shown in \cite[Theorem 1.16]{BEGZ10} that two cohomologous positive closed currents 
$T \sim S$ satisfy
$$
\langle T^n \rangle \leq \langle S^n \rangle 
$$
if $T$ is more singular than $S$ and both currents have {\it small unbounded locus}, i.e. their local potentials are locally bounded outside a closed complete pluripolar set. 
This latter condition seems purely technical and the question makes sense in full generality:

 \begin{quest} [Boucksom-Eyssidieux-Guedj-Zeriahi]
 Assume $T \sim S$ with $T$ more singular than $S$. Is it  true that
$$
\langle T^n \rangle \leq \langle S^n \rangle  \; \; ?
$$
\end{quest}

The construction of the non-pluripolar product makes sense when the cohomology class of $T,S$ is  not necessarily  K\"ahler, but it uses the fact that the manifold is compact and K\"ahler (or at least
belongs to the Fujiki class). It would be interesting to understand whether non-pluripolar products are well-defined on arbitrary compact complex manifolds.

In the sequel we set
$$
MA(\f):=V^{-1} \langle (\omega+dd^c \f)^n \rangle,
\text{ with } 
V=\int_X \langle \omega^n \rangle.
$$
This is a positive Radon measure of total mass at most one.

\begin{defi}
We let ${\mathcal E}(X,\omega)$ denote the set of those functions $\f \in PSH(X,\omega)$ such that
the measure $MA(\f)$ is a probability measure.
\end{defi}

The "finite energy class" ${\mathcal E}(X,\omega)$   has been introduced in \cite{GZ07} and further studied in \cite{CGZ08,BEGZ10,BBGZ13}. Clearly $PSH(X,\omega) \cap L^{\infty}(X) \subset  {\mathcal E}(X,\omega)$, but the finite energy class also contains many unbounded functions: for example if $\f$ is an arbitrary
$\omega$-psh function then $-(-\f-\sup_X \f-1)^\a \in {\mathcal E}(X,\omega)$ for any $0 \leq \a <1$.

The range of the complex Monge-Amp\`ere operator is well understood in this context:

\begin{thm} \label{thm:gz07} \cite{GZ07}
  A probability measure  belongs to the range of the complex Monge-Amp\`ere operator acting
on  ${\mathcal E}(X,\omega)$
if and only if it does not charge pluripolar sets.
\end{thm}

A subtle point is to understand precisely when a function $\f \in DMA(X,\omega)$ belongs to 
${\mathcal E}(X,\omega)$. This has been investigated in \cite{CGZ08} with poor success:

 \begin{quest} [Coman-Guedj-Zeriahi]
 Assume $\f \in DMA(X,\omega)$ is such that $MA(\f)$ does not charge $(\f=-\infty)$.
Does it imply that $\f$ belongs to ${\mathcal E}(X,\omega)$ ?
\end{quest}

The answer is positive when $n=1$. The problem seems more subtle when 
$n \geq 2$, it is again positive  when $e^{\f}$ is continuous \cite[Proposition 1.6]{CGZ08}.

\subsection{Pluripolar measures}

Thanks to the works of Cegrell  and the authors, there is a fairly good understanding of the non-pluripolar  range of the complex Monge-Amp\`ere operator both in domains
of $\C^n$ \cite{Ceg98,Ceg04} and on compact K\"ahler manifolds \cite{GZ07,Din09}.

Despite recent results by Cegrell and his coauthors \cite{ACCH09,ACH15}, the corresponding problem
for measures which do charge pluripolar sets is fully open.
We mention the following particular issue:

 \begin{quest} [Phong]
 Let $X \stackrel{\pi}{\rightarrow} \C\PP^2$ be the blow up of the complex projective plane at one point.
This is a DelPezzo surface, hence its first Chern class is positive.
Pick $\omega \in c_1(X)$  a K\"ahler form,
$D$ a smooth anticanonical divisor, and let $\mu_D=\omega \wedge [D] /V$ denote a normalized  
Lebesgue measure along $D$. Can one find $\f_D \in DMA(X,\omega)$ such that
$$
MA(\f_D)=\mu_D \; \; ?
$$
\end{quest}

The measure $\mu_D$ is pluripolar as it is supported on $D$.
The question is motivated by the fact that $X$ does not admit any K\"ahler-Einstein metric, hence
the K\"ahler-Ricci flow diverges on $X$.
Such a Monge-Amp\`ere potential $\f_D$ could be a cluster point of renormalized
K\"ahler potentials evolving along the K\"ahler-Ricci flow on $X$.

\subsection{The subsolution problems}

We now focus on the case of {\it bounded} and {\it continuous} plurisubharmonic functions.
Although it is the setting of the original breakthrough of Bedford and Taylor \cite{BT76,BT82},
the range of the complex Monge-Amp\`ere operator acting on bounded/continuous  (quasi-)psh functions
is yet not characterized (despite rather sharp sufficient conditions \cite{Kol98}).
We recall the following question from \cite[p. 70]{Bed93}:

 \begin{quest} [Bedford]
 Is the range of the complex Monge-Amp\`ere operator acting on $PSH \cap L_{loc}^{\infty}$
characterized by a uniform integrability condition ?
\end{quest}

When $n=1$, the answer is positive: a probability measure $\mu$  can be written
$\mu=\Delta \f$ for some {\it locally bounded} subharmonic function $\f$ if and only if
it satisfies
$$
SH(\Omega) \subset L^1_{loc}(\mu),
$$
as the reader will check. The integrability condition $PSH(\Omega) \subset L^1_{loc}(\mu)$ is no longer sufficient if $n \geq 2$.

In a bounded hyperconvex domain $\Omega$ of $\C^n$, 
a probability measure $\mu$  can be written
$$
\mu=MA(\f)
\; \; \text{ with } \; \;
 \f_{|\partial \Omega} \equiv 0
$$
if and only if $\mu$ is dominated by the Monge-Amp\`ere measure of a bounded
plurisubharmonic function. This follows from the comparison principle
and the homogeneity of $MA$ in the local context,
$MA(\lambda \f)=\lambda^n \MA(\f)$ \cite{Kol95}. In other words: if there is a bounded subsolution, then there is a bounded solution to the corresponding Dirichlet problem.

 \begin{quest} [Ko\l{}odziej]
 Assume there is a continuous subsolution to the Dirichlet problem, i.e. there is  a continuous
psh function $u$ in $\Omega$ such that $u_{|\partial \Omega} =0$ and
$\mu \leq (dd^c u)^n$ in $\Omega$. Does it imply that there is a continuous solution, i.e.
a continuous psh function $\f$ such that $\f_{|\partial \Omega} =0$ and
$\mu = (dd^c \f)^n$ in $\Omega$ ? 
\end{quest}

The answer is positive in dimension $n=1$ as follows from the following characterization: a probability measure $\mu$  can be written
$\mu=\Delta \f$ for some {\it continuous} subharmonic function $\f$ if and only if
the inclusion
$$
SH(\Omega) \hookrightarrow L^1_{loc}(\mu)
$$
is continuous,
as the reader will check. 

The situation is even more delicate on compact K\"ahler manifolds:

 \begin{quest} [Ko\l{}odziej]
 Let $(X,\omega)$ be a compact K\"ahler manifold. Let $\mu$ be a probability measure which is
locally dominated by the Monge-Amp\`ere measures of bounded plurisubharmonic functions.
Can one find $\f \in PSH(X,\omega) \cap L^{\infty}(X)$ such that 
$\mu=MA(\f)$ ?
\end{quest}

The homogeneity property does not hold in the compact setting. Ko\l{}odziej has proposed an approach in
\cite{Kol05}, but the proof contains a gap.
One can show that the hypothesis on $\mu$ is equivalent to the existence of 
$u \in PSH(X,\omega) \cap L^{\infty}(X)$ and $C>0$  such that
$$
\mu \leq C \, MA(u).
$$
Note in particular  that $\mu$ does not charge pluripolar sets. It follows therefore from 
\cite{GZ07,Din09} that there exists a unique finite energy function 
$\f \in {\mathcal E}(X,\omega)$ such that $\mu=MA(\f)$, and the problem is to show that
$\f$ is bounded below.

A positive answer to this question would imply  that 
the set $MA(PSH(X,\omega) \cap L^{\infty}(X))$ is convex since,
for $u,v \in PSH(X,\omega) \cap L^{\infty}(X)$, 
$$
MA(u)+MA(v) \leq 2^n \, MA \left( \frac{u+v}{2} \right).
$$

\subsection{H\"older continuous solutions}

Let $(X,\omega)$ be a compact K\"ahler manifold of dimension $n$. We let
$$
{\mathcal Hol}(X,\omega):=PSH(X,\omega) \cap \rm{Holder}(X)
$$
denote the set of H\"older continuous $\omega$-psh functions on $X$.

Recall that the Lelong numbers of $\omega$-psh functions are uniformly bounded
by a constant that only depends on $X,n,\{\omega\}$.
It follows therefore from Skoda's integrability theorem \cite{Sko72} (see \cite{Zer01} for the uniform version)
that there exists $\e>0$ such that
$$
\exp(-\e PSH(X,\omega) ) \subset L^1(\mu),
$$
for any smooth volume form $\mu$.

  Dinh-Nguyen-Sibony \cite{DNS10} have shown that 
the same integrability property holds true for any Monge-Amp\`ere measure
$\mu=MA(\f)$ of a H\"older continuous $\omega$-psh function $\f$.
We wonder whether this property actually characterizes $MA({\mathcal Hol}(X,\omega))$:

 \begin{quest} [DDGHKZ]
 Are the following properties equivalent ?

\smallskip

(i) There exists $\e_{\mu}>0$ such that $\exp(-\e_\mu PSH(X,\omega) ) \subset L^1(\mu)$.

\smallskip

(ii) There exists $\f \in {\mathcal Hol}(X,\omega)$ sucht that $\mu=MA(\f)$.
\end{quest}

These properties are indeed equivalent when $n=1$, or when the singularities have toric symmetries
(see \cite{DDGHKZ14}). 

A (technical) characterization of $MA({\mathcal Hol}(X,\omega))$ is established in
\cite{DDGHKZ14}, which solves in particular the subsolution problem in the H\"older case:
a measure belongs to $MA({\mathcal Hol}(X,\omega))$ if and only if it admits a H\"older continuous
subsolution. This property is therefore local. The local H\"older subsolution problem is however
open:

 \begin{quest} [Zeriahi]
Let $\mu$ be a probability measure in a bounded pseudoconvex domain $\Omega \subset \C^n$.
  Assume there is a H\"older continuous subsolution to the Dirichlet problem, i.e. there is  a 
H\"older continuous psh function $u$ in $\Omega$ such that $u_{|\partial \Omega} =0$ and
$\mu \leq (dd^c u)^n$ in $\Omega$. 
Is there  a H\"older continuous psh function $\f$ such that $\f_{|\partial \Omega} =0$ and
$\mu = (dd^c \f)^n$ in $\Omega$ ? 
\end{quest}

\section{Regularity issues}

\subsection{Uniform a priori estimates}

The foundations of an existence and regularity theory for complex Monge-Amp\`ere 
equations with smooth data were laid by Yau \cite{Yau78} and Caffarelli-Kohn-Nirenberg-Spruck \cite{CKNS85}.

 Yau notably provided in \cite{Yau78} a crucial $L^{\infty}$- a priori estimate in the context of
compact K\"ahler manifolds. His method, based on Moser iteration technique can be generalized to work for right hand side data in $L^p$ with $p>n=dim_{\mathbb C} X$. Later on Cheng and Yau suggested that the Alexandrov-Bakelman-Pucci
estimate can be applied to the complex Monge-Amp\`ere equation. A detailed account of their work
has been provided by Bedford \cite{Bed93} and Cegrell-Persson \cite{CP92}. This approach works for $L^p$- right hand side with $p\geq 2$.
This has been since then extended by Ko\l{}odziej who showed in \cite{Kol96} the following:

 \begin{thm} 
  Let $\B$ be the unit ball of $\C^n$ and fix $p>1$. If $\f$ is a smooth plurisubharmonic function
in $\B$ such that $\f_{| \partial \B} =0$, then
$$
|| \f ||_{L^{\infty}(\B)} \leq C ||f||_{L^p(\B)}^{1/n},
$$
where $\det(\f_{j\overline{k}})=f$ and $C$ only depends on $n,p$.
\end{thm}

It is interesting that no PDE proof of the above fact is known.
 \begin{quest} [Blocki-Ko\l{}odziej]
 Find a proof of the above estimate that relies on PDE techniques.
\end{quest}

\subsection{H\"older continuity}

Let $X$ be a $n$-dimensional compact K\"ahler manifold. Let $\omega$ be a semi-positive closed $(1,1)$-form with $V_{\omega}:=\int_X \omega^n>0$, and fix $0 \leq f \in L^p(dV)$, where
$dV$ is a volume form and $p>1$.

When $\omega$ is K\"ahler it has been shown by Ko\l{}odziej in \cite{Kol08} that the equation
$$
(\omega+dd^c \f)^n=f \, dV
$$
admits a unique normalized $\omega$-psh function $\f$ which is H\"older continuous.

This result has been refined in \cite{EGZ09,Din10, Hiep10,DDGHKZ14} (better exponent, $L^p$-property and control in families). The following local $1$-dimensional example shows that one cannot expect
more than H\"older regularity, and how the exponent is linked to the integrability properties of the density:

\begin{exa}
 The function $\f:z  \longmapsto |z|^{2\a} $ is subharmonic in the unit disk of $\C$ with
$$
dd^c \f=\frac{c}{|z|^{2-2\a}} dV=f dV,
$$
with $f \in L^p$ if and only if $p<1/(1-\a)$.
\end{exa}

If the form $\omega$ is merely semi-positive and big (i.e. $\int_X\omega^n>0$) then solutions to the Monge-Amp\`ere equations with $L^p$- right hand sides are continuous  as shown in \cite{EGZ11} by using viscosity techniques (solutions are then unique up to an additive constant). Further regularity is then largely open. In particular the following is an important problem:

\begin{quest} [Eyssidieux-Guedj-Zeriahi]
Assume $\omega$ is merely semi-positive and big and 
$$
(\omega+dd^c \f)^n=f \, dV,\ f\in L^p,\ p>1.
$$

Is $\f$ H\"older continuous on $X$ ?
\end{quest}

It is known that the solution is  H\"older continuous on the ample locus of $\{\omega\}$
\cite{DDGHKZ14}. The problem is thus to understand the regularity at the boundary of the ample locus.

Such equations show up naturally in K\"ahler geometry, when constructing singular K\"ahler-Einstein metrics
on mildly singular varieties (see \cite{EGZ09}). The problem is to understand the asymptotic behavior of these metrics near the singularities; the function $\f$ above occurs as the potential of such a metric. 

The orbifold setting (quotient singularities) is well understood and shows that one cannot expect more than global H\"older continuity of the potentials near the singularities. One can  hope to be able to give precise asymptotics, but this remains a major open problem so far.

\subsection{Totally real submanifolds}

Let $(X, \omega)$ be a compact $n$- dimensional K\"ahler manifold. Recall that a real submanifold $S \subset X$ is called {\it totally real}
if for every point $x\in S$ the real tangent plane $T_xS$ does not contain complex lines.

Let $\mu_S$ denote the normalized $n$-dimensional Hausdorff measure on a smooth totally real manifold $S \subset X$ of maximal real dimension $n$.
It has been shown by Sadullaev in \cite{Sad} that $S$ is non pluripolar, hence $\mu_S$ is the Monge-Amp\`ere measure 
$\mu_S=(\omega+dd^c \f_S)^n$ of a unique function $\f_S \in {\mathcal E}(X,\omega)$ normalized by $\sup_X \f_S=0$.

\begin{quest}[Berman]
What is the optimal regularity of $\f_S$ ?
\end{quest}

One expects $\f$ to be Lipschitz in this case: the torus $S=S^1 \times \cdots S^1 \subset \C\PP^n$ is such that
$$
\f_S=\max_{0 \leq i \leq n} \log |z_i| -\log ||z||.
$$

The potential $\f_S$ is Lipschitz continuous for any real analytic totally real submanifold $S \subset \C\PP^n$.
It has been shown by Sadullaev and Zeriahi \cite{SadZer} that extremal functions of such manifolds
are Lipschitz-continuous.

\subsection{Smooth currents in big cohomology classes} \label{sec:big}

Let $X$ be a compact manifold in the Fujiki class (i.e. $X$ is bimeromorphic to a compact K\"ahler manifold).
It follows from the work of Demailly 
that $X$ admits one (hence many) 
{\it K\"ahler current}(s), i.e. a positive closed $(1,1)$-current $T$ which dominates a positive hermitian form.

\begin{defi}
A cohomology class $\a \in H^{1,1}(X,\R)$ is big if it can be represented by a K\"ahler current.
The ample locus $Amp(\a)$ is the set of points $x \in X$ such that there exists a K\"ahler current
representing $\a$ which is smooth near $x$. 
\end{defi}

The ample locus is a non empty Zariski open subset, as follows from Demailly's regularization result
\cite{Dem92}. The latter insures that any big cohomology class can be represented by a K\"ahler current with
{\it analytic singularities}:

\begin{defi}
A positive closed current $T$ has analytic singularities if it can be locally written $T=dd^c u$, with
$$
u=\frac{c}{2}\log \left[ \sum_{j=1}^s |f_j|^2 \right]+v,
$$
where $c>0$, $v$ is smooth and the $f_j$'s are holomorphic functions.
\end{defi}

The theory of complex Monge-Amp\`ere equations has been extended to big cohomology classes in \cite{BEGZ10}.  
Theorem \ref{thm:gz07} holds in this context (the boundary of the ample locus is pluripolar hence negligible for non-pluripolar measures): given $\mu$ a non-pluripolar probability measure,
there exists a unique finite energy current $T \in \a$ such that 
$$
\langle T \rangle^n=\mu,
$$
where $\langle T \rangle^n$ denotes the non-pluripolar product.
The regularity theory is more involved:

\begin{quest}
Let $\a \in H^{1,1}(X,\R)$ be a big cohomology class, $\mu$  a smooth volume form and let $T$ be the unique finite energy current in $\a$ such that 
$\langle T \rangle^n=\mu$.
 Is $T$ smooth in $Amp(\a)$ ?
\end{quest}

The answer is positive under the extra assumption that $\a$ is {\it nef} \cite{BEGZ10}.
It is also known to hold when $\a=c_1(K_X)$, i.e. when $X$ is a manifold of general type,
but this relies on important recent progresses in birational geometry \cite{BCHM10,EGZ09}.

\medskip

A similar problem concerns Monge-Amp\`ere equations on quasi-projective manifolds.
In geometrical applications, singularities of the right hand side are often of divisorial type i.e. $f\approx ||\sigma||^a$ for some defining section $\sigma$ of a fixed divisor and $a\in\mathbb R$. It is natural to ask whether the singular metric associated to the solution of the Monge-Amp\`ere equation is in fact smooth off the divisor:

\begin{quest}[DiNezza-Lu]
Let $(X,\omega)$ be a compact $n$-dimensional K\"ahler manifold.
 Assume  $\f \in {\mathcal E}(X,\omega)$ solves the equation $(\omega+dd^c \f)^n=f\omega^n$, with $f\in L^1$ being 
smooth off a divisor $D\subset X$. Is the function $\f$ smooth on $X\setminus D$?
\end{quest}

Since the measure $\mu=f\omega^n$ is non-pluripolar, it is known \cite{GZ07,Din09} that there exists a unique
normalized solution $\f \in {\mathcal E}(X,\omega)$, so the point is to study its regularity in $X \setminus D$.
Under some mild additional conditions the answer is affirmative \cite{DiNL14}. The general case however remains open.

\subsection{Pogorelov interior estimate}

A bounded domain $\Omega \subset \C^n$ is {\it hyperconvex} if there exists a negative plurisubharmonic function
$\rho$ in ${\Omega}$ such that the sublevel sets $\{ \rho<c\}$ are relatively compact in $\Omega$ for any $c<0$.
Note that all smooth and strongly pseudoconvex domains are hyperconvex, but they are plenty of hyperconvex non-smooth (and/or not strongly pseudoconvex) domains.

\begin{quest} [Blocki]
Let $\Omega  \subset \C^n$ be a bounded hyperconvex domain.
 Let $u$ be the unique continuous plurisubharmonic function in 
$\Omega$ such that $u_{| \partial \Omega}=0$ and 
$$
(dd^c u)^n=f dV,
$$
with $f$ smooth and bounded away from zero $f \geq \e >0$. Is $u$ automatically smooth ?
\end{quest}

Note that the main result of \cite{Blo96} implies that  any bounded hyperconvex domain   admits {\it one} exhaustion function $u$ with the above properties. 
The problem is open even if one assumes $f$ is constant.

The answer is positive when $\Omega$ is smooth and strongly pseudoconvex \cite{CKNS85}.
It is however crucial for applications to avoid any a priori regularity assumption on $\partial \Omega$,
so as to be able to apply such a result to the open sublevel sets $\Omega_c=\{u<c\}$.

The answer  is also positive for the {\it real Monge-Amp\`ere equation} \cite{Caf90}, the main ingredient being an interior ${\mathcal C}^2$-estimate of Pogorelov \cite{Pog71}.
A complex version of this estimate is not known.
In a similar direction we mention the following:

\begin{quest} [Blocki-Dinew]
 Let $u$ be a plurisubharmonic function in a smooth strictly pseudoconvex bounded domain 
$\Omega \subset \C^n$. Assume $u \in {\mathcal C}^{1,\a}(\Omega),\a>1-\frac2n$ and 
$$
(dd^c u)^n=f dV,
$$
with $f$ smooth and bounded away from zero $f \geq \e >0$. Is $u$ automatically smooth ?
\end{quest}

Example \ref{Pogorelov} below shows that the condition $\a>1-\frac 2n$ cannot be dropped. 
The answer is affirmative for the {\it real Monge-Amp\`ere equation} \cite{Caf90}. In the complex setting the best result is due to Wang \cite{Wa},
who showed that $u$ is smooth assuming $\Delta u\in L^{\infty}(\Omega)$ instead of $u \in {\mathcal C}^{1,\a}(\Omega),\a>1-\frac2n$.

\subsection{Krylov's result}  
It is an interesting and largely open question to determine
when pluripotential solutions are actually {\it strong} solutions, i.e. when the Monge-Amp\`ere operator can be defined pointwise almost everywhere. A special case of the above situation is when the solution is $\mathcal C^{1,1}$ smooth. Hence it is interesting to know when such a regularity holds.

\begin{quest} [Demailly]
Fix $\Omega \subset \C^n$  a smooth strictly pseudoconvex bounded domain, and
$\phi \in {\mathcal C}^{1,1}(\partial \Omega)$.
  Let $u$ be the unique  plurisubharmonic function in $\Omega$, continuous up to the boundary,
such that $u=\phi$ on $\partial \Omega$ and
$$
(dd^c u)^n=0 \; \; 
\text{ in } \; \; \Omega.
$$

Show, by using complex analytic methods, that $u$ is ${\mathcal C}^{1,1}$-smooth in $\Omega$.
\end{quest}

This result has been established by Bedford and Taylor when $\Omega=\B$ is the unit ball \cite{BT76}.
The general case has been established by Krylov in a series of papers culminating in \cite{Kry89}, by using probabilistic methods. 

We refer the reader to \cite{Del12} for a comprehensive presentation of Krylov's approach which remains quite involved. Demailly has suggested on several occasions that one should be able to produce a  proof of Krylov's result by using more standard tools from complex analysis.
We refer the reader to \cite{GZ12} for a brief historical account of related questions.

A related problem concerns ${\mathcal C}^{1,1}$ regularity up to the boundary:

\begin{thm}[Krylov]
Assume that  $\Omega \subset \C^n$  is a ${\mathcal C}^{3,1}$ smooth strictly pseudoconvex bounded domain, and
$\phi \in {\mathcal C}^{3,1}(\partial \Omega)$.
  Let $u$ be the unique  plurisubharmonic function in $\Omega$, continuous up to the boundary,
such that $u=\phi$ on $\partial \Omega$ and
$$
(dd^c u)^n=f \; \; 
\text{ in } \; \; \Omega
$$ 
with $f\geq 0, f^{1/n}\in\mathcal C^{1,1}(\bar{\Omega})$. Then $u\in {\mathcal C}^{1,1}(\bar{\Omega})$.
\end{thm}

Bedford and Fornaess \cite{BeFo} have shown that $u$ need not belong to $\mathcal C^2(\bar{\Omega})$ under the assumptions above. Their construction relies on the boundary behavior of the $\mathcal C^2$ norm.

\begin{quest}
Under what minimial assumptions does $u$ belong to $\mathcal C^2(\Omega)$ ?
\end{quest}

The problem requires to better understand how singularites form {\it in the interior} of $\Omega$.

\subsection{Regularity of the pluricomplex Green function}
Let $\Omega$ be a bounded hyperconvex domain. The pluricomplex Green function with pole at $w\in\Omega$ is defined by
$$g_{\Omega}(z,w)=sup\lbrace u(z)|\ u\in\psh(\Omega),\ u\leq 0,\ u(s)\leq \log ||s-w||+C\rbrace,$$
where the constant $C$ may depend on $u$.

It is easy to prove that for a hyperconvex domain $\Omega$ and a fixed pole $w$ one has
$$\lim_{z\rightarrow\partial\Omega}g_{\Omega}(z,w)=0.$$

\begin{quest}[Blocki]
Is it true that for a bounded hyperconvex domain $g_{\Omega}(z,w)$ converges locally uniformly to zero, as the pole $w$ converges to the boundary of $\Omega$? 
 
\end{quest}

The question is only interesting if $\Omega$ has very rough boundary.

\begin{quest}[Blocki] Let $\Omega$ be smoothly bounded strictly pseudoconvex domain. It is known that $g_{\Omega}\in\mathcal C^{1,1}(\Omega\setminus\lbrace w\rbrace)$. Is there an example where $g$ fails to be in $\mathcal C^2(\Omega\setminus\lbrace w\rbrace)$?
 
\end{quest}

Coman \cite{Co} and independently Edigarian and Zwonek \cite{EdZwo} computed explicitly the pluricomplex Green function {\it with two poles} for the unit ball in $\mathbb C^n$ and the function is not $\mathcal C^2$. No such example with one pole is known. On the other hand Bedford and Demailly \cite{BeDe} have constructed an example which need not be $\mathcal C^2$ smooth {\it up to the boundary} but the interior regularity is open.

\section{Hermitian setting}

In this section we work on a compact complex manifold $X$ that is not necessarily K\"ahler.
The  K\"ahler condition imposes topological and geometrical restrictions.
The reference form $\omega$ is now  the $(1,1)$-form of a {\it hermitian} metric, it is no longer closed.

There always exist such a hermitian metric $\omega>0$, as can be seen by patching together local metrics by partition of unity. The plurisubharmonicity with respect to such a hermitian metric $\omega$ can be defined in an analogous way as in the K\"ahler setting by the inequality $\omega+dd^c\f\geq 0$. Note also that $MA\f)=(\omega+dd^c\f)^n$ makes sense if $\omega$ is merely {\it hermitian} positive, even though it is not true that $\omega$ is locally the $dd^c$ of a plurisubharmonic potential. 

The foundations of pluripotential theory with respect to a general hermitian metric have been developed in \cite{DiKo}.
There are numerous  extra technical difficulties in this hermitian (non K\"ahler) setting. 
The lack of local potentials makes it more difficult to localize arguments in coordinate patches since the outcome is not the usual local Monge-Amp\`ere equation. 
Moreover the Hermitian form $\omega$  being not closed, integrating by parts $(\omega+dd^c\f)^n$ over $X$ 
yields additional terms involving $d\omega$ which have to be handled separately.

\subsection{The total volume}

A substantial difference with the K\"ahler setting is that the total volume $\int_X(\omega+dd^c\f)^n$ depends on $\f$ in general. Such a property is important in various comparison principle arguments \cite{DiKo,KoNg1}.

\begin{quest}[Dinew-Ko\l{}odziej]
Let $(X,\omega)$ be a compact Hermitian manifold of dimension greater or equal to $3$. 
 Characterize the forms $\omega$ for which $\int_X(\omega+dd^c\f)^n=\int_X\omega^n$ for any bounded $\omega$-plurisubharmonic function
 $\f$. 
\end{quest}
 A good candidate for such a condition is the property studied by Guan and Li \cite{GuLi}: namely one assumes that
$$
dd^c\omega=0,\ d\omega\wedge d^c\omega=0.
$$

The same phenomenon leads to the following questions:

\begin{quest}[Tosatti]
Let $(X,\omega)$ be a compact Hermitian manifold. Define $T(X,\omega)$ by
$$
T(X,\omega)=\inf \left\lbrace\int_X(\omega+dd^c\f)^n|\ \f\in\mathcal C^{\infty}(X)\cap PSH(X,\omega) \right\rbrace.
$$
 Is it possible that $T(X,\omega)=0$ for some manifold?
\end{quest}

\begin{quest}[Dinew-Ko\l{}odziej]
Let $(X,\omega)$ be a compact Hermitian manifold. Define $cap(X,\omega)$ by
$$
cap(X,\omega)=\inf \left\lbrace\int_X(\omega+dd^c\f)^n|\ \f\in PSH(X,\omega),\ 0\leq\f\leq 1\right\rbrace.
$$
 Relate the quantity $cap(X,\omega)$ to the geometry of $(X,\omega)$. In particular is this quantity computable in terms of geometric
 data? How does it behave under deformations?
\end{quest}

\subsection{Monge-Amp\`ere equations}
Just like in the K\"ahler case one can pose the problem of solving  the equation
$$
(\omega+dd^c \f)^n=e^cf\omega^n
$$
for a given nonnegative function $f$. The additional constant factor $e^c$ stems from the lack of invariance of the total volumes.

This equation has been solved in the smooth setting by Tosatti and Weinkove \cite{ToWe}:

\begin{thm}[Tosatti-Weinkove]
Let $(X,\omega)$ be a compact Hermitian manifold. Given any smooth strictly positive function $f$ 
there is a unique constant $c$ and a unique function $\f\in\mathcal C^{\infty}(X)\cap PSH(X,\omega)$ satisfying
$$
sup_X\f=0,\ \ \ (\omega+dd^c \f)^n=e^cf\omega^n.
$$
\end{thm}

In the singular setting the following result is due to Ko\l{}odziej and Nguyen \cite{KoNg1}:

\begin{thm}[Ko\l{}odziej-Nguyen]
Let $(X,\omega)$ be a compact Hermitian manifold. The problem
$$
sup_X\f=0,\ \ \ (\omega+dd^c \f)^n=e^cf\omega^n
$$
admits a unique continuous solution $\f$ and a unique constant $c$ for any nonnegative function $0\neq f\in L^p(X,\omega),\ p>1$.
\end{thm}

The higher regularity of solutions was investigated in \cite{KoNg2}. It is based on the following stability result:

\begin{thm}[Ko\l{}odziej-Nguyen]
Let $0\neq f, g$ be two nonnegative functions, such that $f,g\in L^p(\omega^n)$ for some $p>1$. 
Assume $u$ and $v$ solve the problems
$$
u\in PSH(X,\omega)\cap\mathcal C(X),\ sup_Xu=0,\ (\omega+dd^cu)^n=f\omega^n$$
and 
$$v\in PSH(X,\omega)\cap\mathcal C(X),\ sup_Xv=0,\ (\omega+dd^cv)^n=g\omega^n
$$
(the additonal constants have been incorporated in $f$ and $g$). 

{Assume moreover that $f\geq c_0$} for some $c_0>0$. Then
$$||u-v||_{\infty}\leq C||f-g||_p^{\alpha}$$
for any $0<\alpha<\frac 1{n+1}$ and some constant $C=C(X,\omega,c_0,\alpha,p,||f||_p,||g||_p).$
\end{thm}

\begin{quest}[Ko\l{}odziej-Nguyen]
Can one remove the condition $f\geq c_0 >0$ ?
\end{quest}

This would significantly improve the regularity theory in the hermitian setting.

\subsection{Prescribed singularities}
Just as in the K\"ahler case solutions of sequences of complex Monge-Amp\`ere equations with right hand sides clustering along analytic subvarieties should converge weakly to very singular $\omega$-plurisubharmonic functions (compare \cite{Dem92,  DP04}). 

It is expected that the limits should have analytic singularities and hence can serve as singular weights in various $\bar{\partial}$-type problems. The easiest case is when clustering occurs at discrete set of points \cite{Dem92}. Motivated by the K\"ahler setting Tosatti and Weinkove studied the following problem:

\begin{quest}[Tosatti-Weinkove]
Let $z_1,\cdots,z_k$ be a collection of points on a Hermitian manifold $X$ and let $\tau_1,\cdots,\tau_k$ be a collection of positive numbers satisfying $\sum_{j=1}^k\tau_j^n<\int_X\omega^n$. Construct an $\omega$-plurisubharmonic function, so that $\phi(z)\leq \tau_j \log(||z-z_j||)+O(1)$ near $z_j$.
\end{quest}

Tosatti and Weinkowe were able to apply Demailly's method in complex dimension 2 and (under some assumptions) in dimension 3 \cite{ToWe2}. The general case remains open.

\section{Integrability properties}  

\subsection{Desingularization of quasi-plurisubharmonic functions}  

Let $\f$ be a quasi-psh function on  a compact complex hermitian manifold $(X,\omega)$. 
We can assume without loss of generality that $\omega+dd^c \f \geq 0$ is a positive current.

When $n=1$, the Lelong number $\nu(\f,x)$ of $\f$ at $x$ equals the Dirac mass that
the measure $\omega+dd^c \f$ puts at $x$, thus
$E^+(\f):=\{ x \in X \, | \, \nu(\f,x) >0 \}$ is at most countable and
$$
\omega+dd^c \f=\sum_{x \in E^+(\f)} \nu(\f,x) \delta_x +\mu,
$$
where $\mu$ is a positive Radon measure with no atom. Note that for all $\e>0$, the
set $E_\e(\f):=\{ x \in X \, | \, \nu(\f,x) \geq \e \}$ is finite and $\omega+dd^c \f$
can be equally decomposed as
$$
\omega+dd^c \f=\sum_{x \in E_\e(\f)} \nu(\f,x) \delta_x +\mu_\e,
$$
where the sum is finite and $\mu_\e$ is a positive Radon measure with atoms of size less than $\e$.

It follows from Siu's celebrated result \cite{Siu74} that in any dimension $n \geq 2$,
the sets $E_\e(\f):=\{ x \in X \, | \, \nu(\f,x) \geq \e \}$ are closed analytic subsets.
One can similarly decompose
$$
\omega+dd^c \f=\sum_{j \geq 0} m_j [D_j] +R
$$
where $m_j \in \R^+$, $D_j$ is a divisor and $R$ is a closed positive current of bidegree
$(1,1)$ such that $E_+(R)$ does not contain any divisor. The set $E_+(R)$ can however contain infinitely
many  analytic subsets of codimension $\geq 2$.

When $n=\dim_\C \leq 2$, 
it follows from the work of Blel-Mimouni \cite{BM05}, the second author \cite{G05b} and 
Favre-Jonsson \cite{FJ05} that $\f$ can be further decomposed as follows:
for any $\e>0$, there exists $\pi_\e:X_\e \rightarrow X$ a finite composition of blow-ups
such that
$$
\pi_\e^* (\omega+dd^c \f)=[D_\e]+R_\e,
$$
where $D_\e$ is a  simple normal crossing divisor with coefficients in $\R^+$ and 
$R_\e$ is a positive current with Lelong numbers smaller than $\e$ at all points.

 \begin{quest} [Favre-Guedj-Jonsson]
Can one similarly desingularize quasi-plurisubharmonic functions in any dimension $n \geq 3$ ?
\end{quest}

The methods of proofs so far make essential use of the $2$-dimensional setting.
Such a result, if true,  should have many applications.

\subsection{Uniform integrability bounds}
Let $\Omega$ be a bounded hyperconvex domain. The volume of sublevel sets of plurisubharmonic functions in $\Omega$ is exponentially decreasing. It is important to have more precise asymptotics in the case when the Monge-Amp\`ere operator is well defined for such a function. The following is a fundamental result \cite{ACKHZ} which deals with functions living in the class $\mathcal F(\Omega)$ (we refer to \cite{Ceg98} for its definition):

\begin{thm}[Ahag, Cegrell, Hiep, Ko\l{}odziej, Zeriahi]
There is a uniform constant $c_n$, such that for any $s>0$ and any $\f\in\mathcal F(\Omega)$ one has
$$
{\rm Vol}(\lbrace\f<-s\rbrace)\leq c_n\delta_{\Omega}^{2n}(1+s\lambda^{-1})^{n-1}exp(-2ns\lambda^{-1}),
$$
with $\delta_{\Omega}$ denoting the diameter of $\Omega$ whereas $\lambda$ is defined by $\lambda^n=\int_{\Omega}(dd^c\f)^n$. 
\end{thm}

While giving the optimal exponential integrability the result is still not sharp because of the polynomial term in front. 

\begin{quest}[ACHKZ]
 Is it true that (for a different constant $c_n$)
$$
{\rm Vol}(\lbrace\f<-s\rbrace)\leq c_n\delta_{\Omega}^{2n}exp(-2ns\mu^{-1})?
$$
\end{quest}

As explained by Berman and Berndtsson \cite{BB} this problem is linked to the optimal Moser-Trudinger and Brezis-Merle type inequalities in $\Omega$.

\subsection{Lelong numbers of finite energy currents}

Let $X$ be a compact K\"ahler manifold (or a manifold from the Fujiki class).
Let $\omega$ be a smooth semi-positive $(1,1)$ form which is big, i.e.
such that $\int_X \omega^n >0$, where $n=\dim_\C X$.

 \begin{quest} [Berman-Boucksom-Eyssidieux-Guedj-Zeriahi]
 Assume $\f \in {\mathcal E}(X,\omega)$ is a $\omega$-psh function with finite energy.
Is it true that $\f$ has zero Lelong number at all points ?
\end{quest}

This is shown to be the case when $\omega$ is  a K\"ahler form \cite{GZ07}. The argument 
to show that $\nu(\f,a)=0$ requires
the existence of a function $g_a \in PSH(X,\omega) \cap L^{\infty}_{loc}(X \setminus \{x\})$
such that
$$
g(x)=\e \log \rm{dist}(x,a)+O(1),
$$
 where $\e>0$. When $\omega$ is merely semi-positive and big,  the same argument applies at any point of the ample locus of $\{\omega\}$
 (see section \ref{sec:big} for the definition), so the delicate
point is to check that this remains true at the boundary of the ample locus.
A positive answer is known when the situation corresponds to the desingularization of a finite energy
current on a compact normal space \cite[Theorem 1.1]{BBEGZ}:

\begin{thm}[BBEGZ]
Let $(V,\omega_V)$ be a compact normal space and let $\p \in {\mathcal E}(V,\omega_V)$
be a finite energy $\omega_V$-psh function. Let $\pi:X \rightarrow V$ be any resolution of singularities
and set $\omega=\pi^* \omega$, $\f=\p\circ \pi$. 

Then $\p \in  {\mathcal E}(X,\omega)$ has finite energy
and zero Lelong number at all points.
\end{thm}

The ample locus of $\{ \omega\}$ is here $X \setminus \pi^{-1}(V_{sing})$ and the point is to show that
the Lelong number of $\f$ at $x \in \pi^{-1}(V_{sing})$ is controlled by the slope of
$\p$ at $y=\pi(x)$ (see \cite[Appendix A]{BBEGZ}).

This is of course a very particular setting (although already quite useful for geometric applications), and the general case is largely open.

\subsection{Extension problems}

\subsubsection{Extension from a subvariety}
Let $(X,\omega)$ be a compact K\"ahler manifold and $Y \subset X$ a complex submanifold.
It is natural to wonder whether all $\omega$-plurisubharmonic functions on $Y$ are induced by global $\omega$-plurisubharmonic functions:

\begin{quest} [Coman-Guedj-Zeriahi]
Is it true that
$$
PSH \left(Y,\omega_{|Y} \right)=PSH(X,\omega)_{|Y} \; \;  ?
$$
\end{quest}

The answer is positive when $\omega$ is a Hodge form (i.e. it is the curvature form of an ample line bundle), 
as follows from the the main result of \cite{CGZ13}. A result of Matsumura \cite{Mats13} shows that, conversely, such an extension property forces the associated line bundle to be ample.
An alternative approach to the case of Hodge forms has been proposed by Hisamoto \cite{His12}.
It requires a better understanding of 
the Ohsawa-Takegoshi extension theorem with prescription of jets of high order, a problem of independent interest.

For transcendental classes
the answer is easy and positive for smooth strongly $\omega$-psh functions
(see  \cite[Proposition 2.1]{CGZ13}).  A recent result of Collins and Tosatti \cite{CT14}
establishes such an extension result for functions with analytic singularities.
The general case is however largely open.

\subsubsection{Families}

We consider now a holomorphic family of compact K\"ahler manifolds. 
Let $\pi:{\mathcal X} \rightarrow \D$ be a smooth family of compact K\"ahler manifolds (i.e. $\pi$ is a proper holomorphic submersion)
and let ${\mathcal L} \rightarrow {\mathcal X}$ be a pseudoeffective holomorphic line bundle.
We let $X_t=\pi^{-1}(t)$ denote the fibers of $\pi$ and set $L_t={\mathcal L}_{| X_t}$.

Assume $h_0=e^{-\f_0}$ is a singular hermitian metric of $L_0$
on $X_0$ which is "positive", i.e. such that $dd^c \f_0 \geq 0$ in the weak sense of currents.

\begin{quest} [P\u aun]
 Can one find $h=e^{-\f}$ a hermitian metric of  ${\mathcal L}$ on ${\mathcal X}$ which is a semi-positive extension of $\f_0$, i.e.
such that $dd^c \f \geq 0$ and $\f_{|X_0}=\f_0$ ?
\end{quest}

This question is related to the problem of the invariance of plurigenera, a fundamental result due to Siu in the case of projective families (see \cite{Siu98,Siu02,Tsu02,Pau07}).

A related problem is to obtain an effective version of Demailly's extension theorem \cite[Theorem 2.14]{Dem15}.
We thank M.P\u aun for emphasizing this.

\section{Viscosity problems}  

A  viscosity approach complementary to the pluripotential one has been proposed in
\cite{HL09,EGZ11,Wan12}. It raises a number of interesting open questions.

\subsection{Differential tests} 

In this section we fix $\Omega$ a bounded pseudoconvex domain (for example the unit ball) of $\C^n$
and study local properties of plurisubharmonic functions from a viscosity point of view.

\begin{defi}
Let $\f$ be an upper-semi continuous function in $\Omega$ and pick $z \in \Omega$.
A  function is called a differential test from above 
at $z$ if it is ${\mathcal C}^2$-smooth in a neighborhood $V$ of $z$, and such that
$$
(\f-q)(z)=\sup_V (\f-q).
$$
\end{defi}

Let $\f$ be a plurisubharmonic function which belongs to the domain of definition $DMA(\Omega)$.
It is {\it strictly plurisubharmonic} in the classical sense if there exists $\e>0$ such that 
$\f-\e ||z||^2$ is still a plurisubharmonic function. It follows in particular that
$$
(dd^c \f)^n \geq \e^n dV(z)>0.
$$

It is natural to wonder whether the converse implication holds true: if $(dd^c \f)^n>0$, does it imply that
$\f$ is strictly plurisubharmonic ? It is obviously true when $\f$ is ${\mathcal C}^2$-smooth.
The answer is also positive in a dynamical context, as shown by Berteloot and Dupont in \cite{BD05}. It is however false in general as the following example shows:

\begin{exa}\label{Pogorelov}
 The function
$$
\f(z,w)=(1+|z|^2) |w| 
$$
is Lipschitz continuous and plurisubharmonic  in  $\C^2$. It satisfies
$$
(dd^c \f)^2=c dV,
$$
where $dV$ is the Lebesgue volume form and $c>0$ is a normalizing constant, while
$$
dd^c \f \wedge \frac{i}{2} dw \wedge d\overline{w}=c' |w| \, dV,
$$
hence $\f$ is not strictly plurisubharmonic along the line $(w=0)$.
\end{exa}

One can nevertheless show that $\f$ is strictly plurisubharmonic  {\it in the viscosity sense}:

\begin{defi}
We say that $\f$ is strictly plurisubharmonic in the viscosity sense if for all $z \in \Omega$ and
for all differential test from above $q$ at $z$, one has $(dd^c q)_+^n(z)>0$.
\end{defi}

The point in the example above is that there is no differential test from above for $\f$ along the line
$(w=0)$. A simpler example  in any dimension  is provided by $\f(z)=||z||$: there is no differential test from above at the origin; the function is strictly plurisubharmonic  in the viscosity sense, but not in the classical sense.

\begin{quest} [Guedj-Zeriahi]
Let $\f$ be a plurisubharmonic function. What can be said about the set of points at which $\f$ does not admit a differential test from above ?
\end{quest}

It is known that this set has empty interior (see \cite[Remark 2.7]{Zer14}). Is it of zero Lebesgue measure (or capacity) ?

\subsection{Supersolutions}
 
We fix here $f \geq 0$ a continuous density and let $dV$ denote the Lebesgue euclidean measure.
For a smooth hermitian $(1,1)$-form $\a$, we let $ \a_+(x)$  denote
$\a(x)$ if it is non-negative, or zero if $\a(x)$ has at least one negative eigenvalue.

\begin{defi}
 A  supersolution of the equation
$
-(dd^c V)^n+fdV=0
$
is a lower semi-continuous function $v:\Omega \rightarrow \R$ such that
for all $x \in \Omega$ and for every differential test $q$ from below at $x$,
$$
-(dd^c q)(x)^n/dV+f(x)\geq 0.
$$

A  subsolution of the equation
$
-(dd^c V)^n+fdV=0
$
is an upper semi-continuous function $u:\Omega \rightarrow \R$ such that
for all $x \in \Omega$ and for every differential test $q$ from above at $x$,
$$
-(dd^c q) (x)^n/dV+f(x)\leq 0.
$$

A viscosity solution of the equation $-(dd^c V)^n+fdV=0$ is a continuous function which is both a 
super- and a sub-solution.
\end{defi}

It follows from \cite[Proposition 1.5]{EGZ11} that a bounded upper semi-continuous function
$u$ satisfies  $-(dd^c u)^n+fdV \leq 0$ in the viscosity sense if and only if it is plurisubharmonic
and  satisfies $-(dd^c u)^n+fdV \leq 0$ in the pluripotential sense of Bedford-Taylor \cite{BT82}.

It follows from  \cite[Remark 2.9]{EGZ14} that a continuous function 
satisfies $-(dd^c u)^n+fdV= 0$ in the viscosity sense if and only if it is plurisubharmonic
and  satisfies $-(dd^c u)^n+fdV =0$ in the pluripotential sense of Bedford-Taylor.

 The situation is non symmetric and more delicate to interpret for supersolutions.
Observe that a viscosity supersolution is not necessarily plurisubharmonic.
We therefore consider its plurisubharmonic envelope,
$$
P(v):=\left( 
\sup \{ \p \, | \, \p \in PSH(\Omega) \; \& \; \p \leq v \}
\right)^*.
$$
This is the greatest plurisubharmonic function that lies below $v$.

\begin{quest} [Eyssidieux-Guedj-Zeriahi]
Let $v$ be a  continuous viscosity supersolution of the equation
$$
-(dd^c V)^n+fdV=0.
$$

Is $P(v)$ a pluripotential supersolution of the same equation ?
\end{quest}

The answer is positive when $v$ is ${\mathcal C}^2$-smooth, as follows from \cite[Lemma 4.7.1]{EGZ11}.
While the concept of viscosity subsolution makes perfect sense in pluripotential theory, the dual concept of 
viscosity supersolution is more delicate to interpret.

\subsection{Approximation of quasi-plurisubharmonic functions}

\begin{quest} [Eyssidieux-Guedj-Zeriahi]
Let $\omega$ be a smooth semi-positive and big form on a compact K\"ahler manifold $X$.
Is any function $\f \in PSH(X,\omega)$ the decreasing limit of smooth $\omega$-psh functions ? 
\end{quest}

The answer is positive when $\omega$ is K\"ahler, as follows from the work of Demailly \cite{Dem92}
(see also \cite{BK07} for an elementary proof). It is also positive when
$\omega=\pi^* \omega_V$ is the pull-back of a Hodge form $\omega_V$ on a singular variety $V$
by a resolution $\pi=X \rightarrow V$ \cite{CGZ13} (see also \cite[Corollary 2.5]{EGZ13}).

In general, it follows from \cite{EGZ11} that
one can find a  sequence of {\it continuous} $\omega$-psh functions decreasing to $\f$, but
 the  smooth approximation  is largely open.

The viscosity approach proposed in \cite{EGZ11} has been extended to the case of general big cohomology classes in \cite{EGZ13}. Let $\a \in H^{1,1}(X,\R)$ be a big cohomology class and fix 
$\theta \in \a$ a smooth closed $(1,1)$-form.
The function
$$
V_{\theta}:=\sup \{  v \, | \, v \in PSH(X,\theta) \text{ and } \sup_X v\le 0 \}
$$
is an example of $\theta$-psh function with minimal singularities. It satisfies $\sup_X V_{\theta}=0$. We let
$$
P(\a):=\{ x \in X \, | \, V_{\theta}(x) =-\infty \}
\text{ and }
NB(\a):=\{ x \in X\, | \, V_{\theta} \notin L_{loc}^{\infty}(\{x\}) \}
$$
denote respectively the {\it polar locus} and the {\it non bounded locus} of $\a$. 
The definitions clearly do not depend on the
choice of $\theta$ and coincide with the polar (resp. non bounded locus) of any $\theta$-psh function with minimal singularities. Note that $P(\a)$  may a priori not be closed, while $NB(\a)$ is always closed.

\begin{quest} [Eyssidieux-Guedj-Zeriahi]
 Does $P(\a)$ always coincide with $NB(\a)$ ?
\end{quest}

This condition is both necessary and sufficient in order to approximate any $\theta$-psh function $\f$
by a decreasing sequence of {\it exponentially continuous} $\theta$-psh functions $\f_j$,
i.e. such that $e^{\f_j}$ is continuous. While it is easy to construct examples of 
quasi-psh functions $\p$ such that $P(\p)$ is strictly smaller than $NB(\p)$, we do not know a single example of a big cohomology class $\a$ for which $P(\a)$ is smaller than $NB(\a)$.

\subsection{The K\"ahler-Ricci flow}  

The Ricci flow introduced by R.Hamilton is an evolution equation of Riemannian metrics
$(g_t)_{t>0}$ on a fixed riemannian manifold $X$,
 $$
\frac{\partial g_t}{\partial t}=-\Ric(g_t),
$$
starting from an initial riemannian metric $g_0$. When $(X,g_0)$ is actually a K\"ahler manifold,
the K\"ahler condition is preserved along the flow which is thus called the K\"ahler-Ricci flow.

\subsubsection{Smoothing properties: the compact setting}

In this section we fix $X$ a compact K\"ahler manifold and $\a_0 \in H^{1,1}(X,\R)$ a K\"ahler cohomology class. 
Pick $S_0$ a positive closed current representing $\a_0$.
The K\"ahler-Ricci flow starting at $S_0$ is the evolution equation
$$
\frac{\partial \omega_t}{\partial t}=-\Ric(\omega_t),
$$
where $(\omega_t)_{0<<T_{max}}$ is a family of K\"ahler forms and $\omega_t$ weakly converges to $S_0$
as $t \rightarrow 0$.

It is classical \cite{Cao85} that when $S_0$ is a K\"ahler form, such a flow admits a unique solution on a maximal interval of time $[0,T_{max}[$, where
$$
T_{max}:=\sup \left\{ t \geq 0 \, | \, t K_X+\alpha_0   \text{ is K\"ahler}\right\}.
$$
This result has been extended by many authors (see \cite{GZ13} and references therein). 
The best result so far is due to DiNezza and Lu who showed in \cite{DiNL14b} that there exists
a unique such family of K\"ahler forms if and only if $S_0$ has zero Lelong numbers at all points.

When $S_0$ has some non-zero Lelong numbers, DiNezza and Lu construct a maximal family of positive currents
$(\omega_t)_{0<t<T_{max}}$ which are honest K\"ahler forms on a time-increasing family of Zariski open subsets $\Omega_t$ where they satisfy the K\"ahler-Ricci flow equation.
By analogy with Demailly regularization results \cite{Dem92,Dem93}, it is natural to expect that this family
has analytic singularities.

\begin{quest} [DiNezza-Lu]
Does the K\"ahler-Ricci flow produce approximants with analytic singularities ?
\end{quest}

\subsubsection{Monge-Amp\`ere flows in pseudoconvex domains}

The study of complex Monge-Amp\`ere flows in pseudoconvex domains 
has been initiated in \cite{HLi10}. A theory of weak solutions has been proposed in \cite{EGZ14}, using a viscosity approach,
and the smoothing properties have been investigated by Do in \cite{Do15} 
who showed in particular the following:

\begin{thm} [Do]
 Let $\Omega$ be a bounded smooth strictly pseudoconvex domain of $\C^n$
 and $T>0$. Let $u_0$ be a bounded plurisubharmonic function defined on a neighborhood 
$\tilde{\Omega}$ of $\overline{\Omega}$.
 Assume that $\f \in {\mathcal C}^{\infty} (\overline{\Omega} \times [0,T))$ 
satisfies $\f(z,t)=u_0(z)$ for $z \in \partial \Omega$. 
Then there exists a unique function $u \in {\mathcal C}^{\infty} (\overline{\Omega} \times (0,T))$ 
 such that 

(i) for all $t \in (0,T)$, $z \mapsto u_t(z)=u(z,t)$  is a strictly plurisubharmonic function in $\Omega$;

(ii) $\dot{u}=\log \det (u_{i\overline{j} })$  on $\Omega \times (0,T)$;

(iii) $u=\f$ on $\partial \Omega \times (0,T)$;
 
(iv) $\lim_{t \rightarrow 0} u(z,t)=u_0(z)$, for all $z \in \overline{\Omega}$.
 
\noindent  Moreover $u$ is uniformly bounded and  
$u_t$ converges to $u_0$ in capacity when $t \rightarrow 0$.
\end{thm}
 
Do actually deals with more general flows which contain, as a particular case, the local version of the K\"ahler-Ricci flow. These flows can thus be run from a bounded plurisubharmonic initial data.
The following question is natural, in view of the corresponding global results from \cite{GZ13,DiNL14b}:

\begin{quest} [Guedj]
Can one run the above complex Monge-Amp\`ere flow from an arbitrary plurisubharmonic function ?
\end{quest}

\subsubsection{Viscosity comparison principles}  

Let $X$ be a compact K\"ahler space.
The K\"ahler-Ricci flow on $X$ can be reduced to a parabolic scalar equation.
When $X$ is a compact mildly singular variety, this equation is 
a complex degenerate Monge-Amp\`ere flow, as studied for example in \cite{ST09,GZ13,EGZ14,EGZ14b}.
We consider here the  flows
$$
 (\omega_t +dd^c \f_t)^n =  e^{\dot{\f}_t + F (t,x,\f) } \mu (t,x) 
$$
on $[0,T] \times X$, where
  \begin{itemize}
   \item $T \in ]0,+\infty]$;
  \item $\omega=\omega(t,x)$ is a continuous family of semi-positive and big $(1,1)$-forms on $X$, 
  \item $F (t,z,r)$ is continuous in $[0,T[ \times \Omega \times \R$ and non decreasing in $r$,
  \item $\mu (t,z)\geq 0$ is a bounded continuous volume form on $X$,
   \item $\f : X_T :=  [0,T[ \times X \rightarrow \R$ is the unknown function, with $\f_t: = \f (t,\cdot)$. 
   \end{itemize}
We assume that the forms $\omega_t \geq \theta$ are uniformly bounded below by a fixed 
semi-positive and big $(1,1)$-form $\theta$.

A central problem in the viscosity approach is to establish the {\it comparison principle}:
recall that a bounded upper semi-continuous function $u$ is a viscosity subsolution if it
satisfies
$$
 (\omega_t +dd^c \f_t)^n \geq  e^{\dot{\f}_t + F (t,x,\f) } \mu (t,x) 
$$
in the viscosity sense. One defines, by duality, supersolutions as those lower semi-continuous functions
that verify the reverse inequality $\leq$ in the sense of viscosity.

\begin{quest} [Eyssidieux-Guedj-Zeriahi]
Assume $u$ (resp. $v$) is a viscosity subsolution (resp. supersolution) to the above degenerate complex 
Monge-Amp\`ere flow such that
$u \leq v$ on the parabolic boundary $\{0\} \times X$.
Is it true that that $u \leq v$ on $[0,T] \times X$ ?
\end{quest}

Such a comparison principle has been established in \cite{EGZ14b} in several particular cases of interest.
The above general case however remains  open
(note that a general comparison principle has been established in \cite{EGZ14} in the context of 
pseudoconvex domains of $\C^n$).

\end{document}